\documentclass[11pt,leqno]{article}

\usepackage{amstex,amscd,amssymb,theorem}

\long\def\comment#1\endcomment{}

\makeatletter
\begingroup
\gdef\th@dotted{\normalfont\itshape
  \def\@begintheorem##1##2{%
        \item[\hskip\labelsep \theorem@headerfont ##1\ ##2.]}%
\def\@opargbegintheorem##1##2##3{%
   \item[\hskip\labelsep \theorem@headerfont ##1\ ##2\ (##3).]}}
\endgroup
\makeatother

\theoremstyle{dotted}

\newtheorem{theorem}{Theorem}[section]
\newtheorem{lemma}[theorem]{Lemma}
\newtheorem{corr}[theorem]{Corollary}
\newtheorem{prop}[theorem]{Proposition}
\newtheorem{conjecture}[theorem]{Conjecture}

\makeatletter
\begingroup
\gdef\th@upshape{\normalfont
  \def\@begintheorem##1##2{%
        \item[\hskip\labelsep \theorem@headerfont ##1\ ##2.]}%
\def\@opargbegintheorem##1##2##3{%
   \item[\hskip\labelsep \theorem@headerfont ##1\ ##2\ (##3).]}}
\endgroup
\makeatother

\theoremstyle{upshape}

\newtheorem{defn}[theorem]{Definition}
\newtheorem{remark}[theorem]{Remark}

\makeatletter
\renewcommand{\subsection}{\@startsection{subsection}{2}{0pt}{-3ex
plus -1ex minus -0.2ex}{-2mm plus -0pt minus
-2pt}{\normalfont\bfseries}} \makeatother

\makeatletter
\@addtoreset{equation}{section}
\makeatother

\newcommand{\proof}[1][Proof.]{\smallskip\noindent{\em #1}}
\def\endproof{\hfill\ensuremath{\square}\par\medskip}

\def\eqref#1{\thetag{\ref{#1}}}

\let\latexref=\ref
\def\ref#1{{\normalfont{\latexref{#1}}}}

\newcommand{\wt}{\widetilde}

\newcommand{\idot}{{\:\raisebox{1pt}{\text{\circle*{1.5}}}}}
\newcommand{\highdot}{{\:\raisebox{3pt}{\text{\circle*{1.5}}}}}
\newcommand{\C}{{\mathbb C}}
\newcommand{\Z}{{\mathbb Z}}
\newcommand{\Q}{{\mathbb Q}}

\newcommand{\K}{{\cal K}}

\newcommand{\m}{{\mathfrak m}}

\newcommand{\Spec}{\operatorname{Spec}}

\newcommand{\cchar}{\operatorname{\sf char}} 
\renewcommand{\dim}{\operatorname{\sf dim}} 
\newcommand{\codim}{\operatorname{\sf codim}} 
\newcommand{\cl}{\operatorname{\sf cl}}
\newcommand{\tr}{\operatorname{\sf tr}}
\newcommand{\pt}{\operatorname{\sf pt}}
\newcommand{\disc}{\operatorname{\sf disc}}
\newcommand{\age}{\operatorname{\sf age}}

\newcommand{\Stab}{\operatorname{\sf Stab}}
\newcommand{\Cent}{\operatorname{\sf Cent}}

\title{McKay correspondence for symplectic quotient singularities}

\author{D. Kaledin}

\begin{document}

\maketitle

\tableofcontents

\section*{Introduction}

Let $V$ be a finite-dimensional complex vector space, and let $G
\subset SL(V)$ be a finite subgroup. The quotient variety $X = V/G$
is usually singular. If $\dim V = 2$, then there exists a canonical
smooth projective resolution $\pi:Y \to X$ of the singular variety
$X$. This resolution is crepant, which in this case means that the
manifold $Y$ has trivial canonical bundle (we recall the precise
definition in Definition~\ref{crepant}). The fiber $\pi^{-1}(0)
\subset Y$ over the singular point $0 \in X = V/G$ is a rational
curve, whose components are numbered by the non-trivial conjugacy
classes in the finite group $G$. The homology classes of these
components freely generate the homology group $H_2(Y,\Q)$. All the
other homology groups of the manifold $Y$ are trivial (except for
$H_0(Y,\Q) \cong \Q$). 

This situation has been described by J. McKay in \cite{M}. It is
known as {\em the McKay correspondence}.

In higher dimensions the picture splits into two parts. One can
consider two separate questions.
\begin{enumerate}
\item When does the quotient $X = V/G$ admit a smooth crepant
resolution $Y \to X$?
\item Assuming there exists a smooth crepant resolution $Y \to X$,
what can one say about the homology space $H_\idot(Y,\Q)$?
\end{enumerate}

In the case $\dim V = 3$, the first question has been solved
independently by several people (for the references see
\cite{IR}). In this case no additional conditions are
necessary. Every quotient $X = V/G$ of the $3$-dimensional vector
space $V = \C^3$ by a finite subgroup $G \subset SL(V)$ admits a
smooth crepant resolution $Y \to X$.

The second question has been considered by Y. Ito and M. Reid in
their paper \cite{IR}. For $\dim V = 3$, they proved the same result
that holds for $\dim V = 2$: the total homology space
$H_\idot(Y,\Q)$ of any smooth crepant resolution $Y \to X$ of the
quotient $X = V/G$ has a natural basis numbered by the conjugacy
classes in the group $G \subset SL(V)$. They called this {\em the
McKay correspondence in dimension $3$}.

Roughly speaking, the proof goes as follows. There are two
non-trivial homology groups, namely, $H_2(Y, \Q)$ and $H_4(Y,
\Q)$. Ito and Reid consider them separately. The picture for
$H_2(Y, \Q)$ is similar to the classical case $\dim V = 2$: the
group $H_2(Y, \Q)$ has a natural basis given by the fundamental
classes of the components of the exceptional divisor $E \subset Y$.
Ito and Reid use the theory of valuations to identify the set of
these components with a subset of the conjugacy classes in $G$ (the
set of so-called {\em junior} conjugacy classes). Then they apply
Poincare duality, and identify the rest of the conjugacy classes
with elements in a basis of $H_4(Y, \Q)$.

Additionally, \cite{IR} contains some discussion of the general case
$\dim V > 3$. In particular, the construction for $H_2(Y, \Q)$
generalizes straightforwardly. Unfortunately, in higher dimensions
this is not enough to describe the total homology space $H_\idot(Y,
\Q)$.

The discussion in \cite{IR} has been clarified and further
generalized by Reid in \cite{R} (we also refer the reader to this
paper for a comprehensive bibliography and the history of the
question, which dates back both to \cite{M} and to the physical
paper \cite{vafa}). Reid considered the quotient $V/G$ of an
arbitrary vector space $V$ by a finite subgroup $G \subset
SL(V)$. He assumed that $V/G$ admits a smooth crepant resolution $Y
\to X$. In this situation, he proposed two versions of the
generalized McKay correspondence -- a homological and a
cohomological one. He also suggested a hypothetical construction for
each of these correspondences.  Philosophically, the main
``homological'' claim is the same as for $\dim V = 2$ and $\dim V =
3$.

\begin{conjecture}
Let $X = V/G$ be quotient of a complex vector space $V$ by a finite
subgroup $G \subset SL(V)$. Assume given a smooth crepant resolution
$Y \to X$.

Then the homology space $H_\idot(Y,\Q)$ admits a ``natural'' basis
numbered by conjugacy classes of elements $g \in G$.
\end{conjecture}

As Reid himself put it, his conjecture should be regarded as ``a
pointer towards the truth rather than the truth itself''. In
particular, the precise meaning of the word ``natural'' is not
known. The only known compatibility condition the good McKay
correspondence must satisfy is that the degree of the homology class
associated to the conjugacy class of an element $g \in G \subset
SL(V)$ must coincide with the so-called {\em age} of the element $g$
(age is a canonical integer associated to every automorphism $g \in
SL(V)$ of a vector space $V$; we recall the definition in
Subsection~\ref{cyclic}).

There is a lot of experimental evidence which supports this
conjecture, and there are some specific pairs $\langle V, G \rangle$
for which such a natural basis has been constructed. Moreover, one
can check the validity of the conjecture by suspending the quest for
a basis in the homology, and asking instead, whether the {\em
dimensions} of the homology groups are correct. This has been proved
recently by V. Batyrev (\cite{Ba}). Another proof, based on ideas of
M. Kontsevich, has been given by J. Denef and F. Loeser (\cite{DL}).

However, if one asks for an actual basis in the homology groups,
then our knowledge is much more limited. Aside from the case $\dim V
= 3$ solved in \cite{IR}, the only known situation where there is
some sort of a general approach is the case of the abelian group
$G$. Here one can use the methods of toric geometry with some
success. In particular, the dimensions of homology groups have been
known for a long time (\cite{BD}). There has also been much progress
on the question of the existence of smooth crepant resolutions; for
this we refer the reader to the paper \cite{DHZ}.

\bigskip

In this paper, we impose a different assumption on the pair $\langle
V,G \rangle$. Namely, we assume that the vector space $V$ is
equipped with a non\--de\-ge\-ne\-rate symplectic form, and that the
finite subgroup $G \subset Sp(V) \subset SL(V)$ preserves not only
the volume form in $V$, but also the symplectic form.

\bigskip

This is suggested by the analogy between quotient singularities
$V/G$ by finite subgroups $G \subset SL(V)$ and the compact
Calabi-Yau manifolds, that is, complex manifolds $X$ with trivial
canonical bundle $K_X$. The quotients $V/G$ by subgroups $G \subset
Sp(V)$ are analogous to compact complex manifolds which are not only
Calabi-Yau, but also hyperk\"ahler.  Just as the theory of
hyperk\"ahler manifolds is in many respects different from the
general theory of Calabi-Yau manifolds, one expects that symplectic
quotient singularities $V/G$ have many interesting special
properties and deserve a closer study. Such a study has been
proposed and initiated recently by A. Beauville in his paper
\cite{B}.

Independently, in a recent paper \cite{V} M. Verbitsky has studied
which of the symplectic quotient singularities $X = V/G$ admit a
smooth crepant resolution. He obtained a strong necessary condition
on the subgroup $G \subset Sp(V)$.

\def\refsemismall{\cite[Proposition~4.4]{K}}

\bigskip

Here we consider the symplectic quotient singularities from the
point of view of the homological McKay correspondence. One expects
significant simplifications due to the following fact: every smooth
crepant resolution $\pi:Y \to X$ of a symplectic quotient
singularity $X = V/G$ is semismall (\refsemismall). This means that
for every irreducible algebraic cycle $Z \subset Y$ we have
$$
\codim \pi(Z) \leq 2 \codim Z.
$$
Using this fact and some simple features of the geometry of the
quotient $X = V/G$, one can show that
\begin{itemize}
\item there is only a finite number of irreducible cycles $Z \subset
Y$ for which we have the equality, $\codim \pi(Z) = 2\codim Z$, and
\item the fundamental classes of these ``maximal'' cycles provide a
natural basis in the total homology space $H_\idot(Y, \Q)$.
\end{itemize}
One hopes then to apply the valuation-theoretic methods of \cite{IR}
to these maximal cycles, and obtain the homological McKay
correspondence not only for $H_2(Y, \Q)$, but for the total homology
space $H_\idot(Y, \Q)$.

We have found out that this program indeed works, and in a very
straightforward way. Thus, for symplectic vector spaces, the
homological McKay correspondence proposed in \cite{R} is not just a
pointer to the truth, but the whole truth. The construction
described by Reid in \cite{R} works beautifully and gives the
desired basis in the homology space $H_\idot(Y,\Q)$. This is the
subject of the present paper.

\bigskip

\noindent
The paper is organized in the following way. 
\begin{itemize}

\item The key element in the Ito-Reid construction is the technique
of valuations. This seems to be a standard approach in modern
birational geometry, but we haven't been able to find a convenient
reference. Since all the facts we need are really quite simple, we
have given a brief general overview of valuations in Section~1. For
the complete proofs see e.g. \cite{bour}. There is also a good
overview in \cite{IR}.

\item In Section~2 we introduce the vector space $V$ and the finite
group $G$ and describe, following Reid \cite{R}, \cite{YPG}, the
so-called {\em monomial valuations} of the quotient variety $X =
V/G$. At this point we are able to state our main result in the
precise form (Conjecture~\ref{main}).

\item In Section~3 we recall certain facts from \cite{K} on the
geometry of an arbitrary smooth crepant resolution $Y \to X$. After
that, we introduce {\em maximal cycles} in the resolution $Y$, and
we prove that maximal cycles in $Y$ are in one-to-one correspondence
with the conjugacy classes in the finite group $G$.

\item Finally, in the purely topological Section~4 we prove that
holomogy classes of maximal cycles in $Y$ form a basis in the
homology space $H_\idot(Y,\Q)$.
\end{itemize}

\subsection*{A rough outline of the proof.} We would like now to
give a rough outline of the proof of our main theorem, freely using
the terminology of the paper \cite{R}. This is completely
independent from the rest of the paper. The reader unfamiliar with
\cite{R} is advised to skip the rest of the introduction and turn to
the main body of the paper, where we give all the necessary
definitions.

First of all, everything in the paper works over an arbitrary
algebraically closed field $k$ of characteristic $\cchar k =
0$. However, as noted in \cite{IR}, the natural formulation of the
homological McKay correspondence would replace conjugacy classes of
elements $g \in G$ with conjugacy classes of maps $\mu_r \to G$,
where $\mu_r \subset K^*$ is the cyclic group of roots of unity of
order $r$ ($r$ runs over all orders of elements in $G$). We have
adopted this point of view, which indeed makes some proofs easier,
and some formulas become more natural. The reader who is mostly
interested in varieties over the complex numbers, $k = \C$, may
simply replace everywhere a homomorphism $g:\mu_r \to G$ with the
image $g(\varepsilon) \in G$ of the canonical generator $\varepsilon
= \exp\left(\frac{2\pi\sqrt{-1}}{r}\right) \in \mu_r$.

Now, the homological McKay correspondence that we prove works in the
following way.
\begin{enumerate}
\renewcommand{\labelenumi}{(\Alph{enumi})}

\item Conjugacy classes of maps $g:\mu_r \to G$ \qquad
$\Longleftrightarrow$

\item Valuations of the function field $K(X) = K(Y)$ \qquad
$\Longleftrightarrow$

\item Maximal cycles $Z \subset Y$ in the resolution $Y \to X = V/G$
\qquad $\Longleftrightarrow$

\item A basis in the homology space $H_\idot(Y,\Q)$.
\end{enumerate}
The precise definition of a maximal cycle is
Definition~\ref{max.cycle}. As we have noted earlier, due to the
semi-smallness of the resolution $Y \to X$, the correspondence
between \thetag{C} and \thetag{D} becomes a simple topological fact
from the theory of perverse sheaves. We explain this in Section~4.

To pass from \thetag{A} to \thetag{B}, one takes, as suggested in
\cite{R}, the so-called {\em monomial valuation} $v_g$ associated to
a (conjugacy class of) a map $g:\mu_r \to G$. This construction is
probably well-known, although we haven't been able to find a
convenient reference, and describe it from scratch in Section~2. In
the process, we prove a certain technical property of monomial
valuations, Lemma~\ref{mon.min}. This lemma is quite simple, but
very important for the rest of the paper. We do not know whether it
is a new result.

To pass from \thetag{B} to \thetag{C}, one takes the center
$\delta(v) \in Y$ of the monomial valuation $v = v_g$, and the
Zariski closure $Z_g = \overline{\delta(v)} \subset Y$. Again, this
is suggested in \cite{R}. Reid indicates that in the general case
this is not enough, and one also has to do some ``mysterious
cookery''. However, none is needed in the symplectic case $G \subset
Sp(V)$. Moreover, the monomial valuation $v_g$ simply coincides with
the $Z_g$-adic valuation of the function field $K(Y) = K(X)$
(Proposition~\ref{mon=>max}).

The correspondence $\thetag{C} \Rightarrow \thetag{B} \Rightarrow
\thetag{A}$ uses the ramification theory. This is the technique used
in \cite{IR} (but for some reason it is not even mentioned in the
overview given in \cite{R}). To pass from a cycle $Z \subset Y$ to a
map $g:\mu_r \to G$, one takes the $Z$-adic valuation $v = v_Z$ of
the function field $K(Y) = K(X)$ and considers the Galois extension
$K(V)/K(X)$ with Galois group $G$. The number $r$ is the
ramification index of the extension $K(V)/K(X)$ in the valuation
$v$, and the cyclic subgroup $g(\mu_r) \subset G$ is the associated
inertia subgroup.

This construction is the subject of Section~3. The main result is
Proposition~\ref{max=>mon}, which claims that the correspondence
$\thetag{C} \Rightarrow \thetag{B} \Rightarrow \thetag{A}$ is indeed
inverse to the correspondence $\thetag{A} \Rightarrow \thetag{B}
\Rightarrow \thetag{C}$ constructed by Reid.

\renewcommand{\labelenumi}{{\normalfont(\roman{enumi})}}

\section{Valuations: a brief overview}

\subsection{Generalities.}

Fix once and for all a base field $k$ which is algebraically closed
and of characteristic $\cchar k = 0$. Throughout the paper, we will
only consider fields $K$ that contain $k$ and such that the
extension $K/k$ is of finite degree of transcendence $\deg\tr K/k <
\infty$.

Let $K$ be such an extension. By a (discrete ultrametric) {\em
valuation} $v$ of the field $K$ we understand a surjective map
$v:K^* \to \Z$ which satisfies
\begin{align}
v(xy) &= v(x) + v(y), \tag{V1}\\
v(x+y) &\geq \min(v(x),v(y)).\tag{V2}
\end{align}
for every $x,y \in K$. Moreover, we assume that $v(x) = 0$ for every
constant $x \in k \subset K$.

Every valuation $v$ of a field $K$ induces a topology on $K$. We
denote by $K_v$ the completion of the field $K$ with respect to this
topology. The {\em ring of integers} $O_v \subset K_v$ associated to
the valuation $v$ consists of $0 \in K$ and all elements $x \in K_v$
such that $v(x) \geq 0$. The ring $O_v$ is a regular local ring of
dimension $1$. The maximal ideal $\m_v \subset O_v$ is generated by
an element $u \in O_v$ with $v(u) = 1$ called the {\em uniformizing
element}. The quotient $k_v = O_v/\m_v$ is a field, it is called the
{\em residue field} of the valuation $v$. The spectrum $\Spec O_v$
is a smooth curve with two points: the generic point $\Spec K_v$ and
the special point $\Spec k_v$.

The module $\Omega^1(O_v/k)$ of K\"ahler differentials is a free
$O_v$-module of rank $n = \deg\tr K/k$, so that one can consider the
modules $\Omega^i(O_v/k) = \Lambda^i(\Omega^1(O_v/k))$ of $i$-forms
on $\Spec O_v$ for every $i \geq 1$. Define a decreasing filtration
$F^\highdot$ on $\Omega^i(O_v/k)$ by
$$
F^p\Omega^i(O_v/k) = 
\begin{cases} 
\Omega^i(O_v/k), &\quad p = 0,\\
\m_v \cdot \Omega^i(O_v/k) + \Omega^{i-1}(O_v/k) \wedge du \subset
\Omega^i(O_v/k), &\quad p = 1,\\
\m_v^{p-1} \cdot F^1\Omega^i(O_v/k), &\quad p \geq 2,
\end{cases}
$$
and set 
\begin{equation}\label{val.forms}
v(\alpha) = \max \left\{p | \alpha \in F^p\Omega^i(O_v/k)\right\}
\end{equation}
for every $i$-form $\alpha \in \Omega^i(O_v/k)$. This extends the
valuation $v:K^* \to \Z$ to a function $v:\Omega^i(O_v/k) \setminus
\{0\} \to Z$ defined for every $i \geq 0$. We have
$$
v(f\alpha) =v(f) + v(\alpha), \quad f \in O_v, \alpha \in
\Omega^\highdot(O_v/k),
$$
so that the function $v$ extends even further to a function
$v:\Omega^\highdot(K_v/k) = \Omega^\highdot(O_v/k) \otimes_{O_v} K_v
\to \Z$. By restriction, we obtain a function
$v:\Omega^\highdot(K/k) \setminus \{0\} \to \Z$. This function
satisfies the inequalities
\begin{align*}
v(\alpha \wedge \beta) &\geq v(\alpha) + v(\beta), \quad \alpha,\beta
\in \Omega^\highdot(O_v/k)\\
v(\alpha + \beta) &\geq \min (v(\alpha), v(\beta) ) \quad \alpha,\beta
\in \Omega^\highdot(O_v/k)\\
v(df) &\geq v(f), \quad f \in O_v,
\end{align*}
where $d:O_v \to \Omega^1(O_v/k)$ is the de Rham differential.

\subsection{Valuations of algebraic varieties.} 

If $X$ is an irreducible algebraic variety, then valuations of $X$
are by definition valuations of the field $K(X)$ of rational
functions on $X$. Thus the set of valuations is a birational
invariant of the variety $X$.

A valuation $v$ of an irreducible algebraic variety $X$ induces a
map 
$$
\eta_v:\Spec K_v \to X. 
$$
If this map extends to a map $\eta_v:\Spec O_v \to X$, then the
point $\delta(v) = \eta_v(\Spec k_v) \in X$ is called the {\em
center} of the valuation $X$, and we say that the valuation $v$ has
a center.

By the valuative criterion of properness, every valuation $v$ of a
proper irreducible algebraic variety $X$ has a center. However, the
following situation will be more important for us. Assume given
irreducible algebraic varieties $X, \wt{X}$, not necessarily proper,
and a proper dominant map $\pi:\wt{X} \to X$ which is generically
one-to-one. Then $\pi$ is a birational isomorphism, so that
valuations of $X$ and $\wt{X}$ are the same. We note that if a
valuation $v$ has a center $\delta(v) \in X$ in $X$, then it also
has a center $\wt{\delta}(v) \in \wt{X}$ in $\wt{X}$, and we have
$\wt{\delta}(v) \in \pi^{-1}(\delta(v)) \subset \wt{X}$. This again
follows from the valuative criterion.

We also note that if the variety $X = \Spec O$ is affine, then a
valuation $v$ of the variety $X$ has a center if and only if $v(f)
\geq 0$ for every algebraic function $f \in O$ on $X$. In this case
the center $\delta(v) \in X$ is the generic point of the closed
subvariety $Z \subset X$ defined by the prime ideal
$$
I_v = \left\{f \in O| v(f) > 0\right\} \subset O.
$$
Assume from now on that the variety $X$ is normal. Then $X$ is
non-singular in codimension $1$. For every irreducible divisor $Z
\subset X$, $\codim Z = 1$, the local ring of functions on $X$ near
the generic point $z \in Z$ of the divisor $Z$ is a discrete
valuation ring. We denote the corresponding valuation of the variety
$X$ by $v_Z$. For any valuation $v$ of $X$ with center at $z$, the
associated map $\Spec O_v \to X$ is an open embedding into the
formal neighborhood of $Z = \overline{z} \subset X$. Therefore $O_v
\cong O_{v_Z}$, and the valuation $v$ must coincide with $v_Z$.

More generally, let $z \in X$ be an arbitrary scheme point in an
irreducible algebraic variety $X$, and let $Z = \overline{z} \subset
X$ be the Zariski closure of the point $x$. Even if $\codim Z \geq
2$, one can still define the valuation $v_Z$ by blowing up $Z
\subset X$ and taking $V_Z$ to be the valuation $v_E$ associated to
the exceptional divisor in the blow-up. The valuation $v_Z$ is
called the {\em $Z$-adic} valuation of the variety $X$.  The
$Z$-adic valuation is obviously centered at $\delta(v_Z) = z \in X$.

\subsection{Discrepancy.}

Assume now that the irreducible algebraic variety $X$ is normal and
has {\em Gorenstein singularities}, that is, the canonical bundle
$K_U$ of the non-singular part $U \subset X$ extends to a line
bundle $K_X$ on the whole $X$, and the extended bundle $K_X$ is
locally trivial. In this situation, there is an important invariant
of centered valuations of $X$ called the {\em discrepancy}. To
define it, assume given a valuation $v$ of the variety $X$ with
center $\delta(v) \in X$ and choose a form $\omega \in
\Omega^n(K/k)$ of top degree at the generic point $\Spec K \subset
X$ which defines a trivialization of the canonical bundle $K_X$ near
the point $\delta(v) \in X$.

\begin{defn}\label{disc.defn}
The {\em discrepancy} of the valuation $v$ in the variety $X$ is the
number
$$
\disc(v,X) = v(\eta_v^*\omega) - 1.
$$
\end{defn}

Here $v$ is the valuation extended to the module $\Omega^n(K_v/k)$
of forms of top degree as in \eqref{val.forms}, and $\eta_v:\Spec
O_v \to X$ is the map associated to the centered valuation $v$, so
that $\eta_v^*\omega$ is an $n$-form $\omega \in \Omega^n(K_v/k)$ on
the generic point $\Spec K_v \subset \Spec O_v$ of the curve $\Spec
O_v$. When there is no danger of confusion, we will write $\disc(v)$
instead of $\disc(v,X)$.

The number $\disc(v)$ does not depend on the choice of
trivialization $\omega$. Indeed, any two trivializations
$\omega_1,\omega_2$ differ by a local function $f$ on $X$ near the
point $\delta(v) \in X$ which does not vanish at $\delta(v) \in
X$. Therefore $v(f) = 0$ and $v(\omega_1) = v(\omega_2)$. A
valuation $v$ of the variety $X$ with $\disc v = 0$ is called {\em
crepant}.

For the valuation $v = v_Z$ associated to a divisor $Z \subset X$,
$\codim Z = 1$, the map $\eta_v$ is an open embedding. Therefore
$\pi^*\omega$ is a generator of the module $\Omega^n(O_v/k)$ of top
degre forms on $\Spec O_v$. Thus every divisorial valuation $v_Z$ is
crepant, $\disc(v_Z) = 0$.

This can be generalized to the situation when the variety $X$ is
smooth near the generic point $z$ of an irreducible subvariety $Z
\subset X$ of arbitrary codimension. The discrepancy of the $Z$-adic
valuation $v_Z$ is easy to compute by the adjunctions formula for
the blowup $\wt{X} \to X$ of $X$ at $Z \subset X$. It equals
\begin{equation}\label{smooth.disc}
\disc(v_Z) = \codim Z - 1.
\end{equation}
In fact, as we shall see in Corollary~\ref{adic.minimal}, this is
the smallest possible discrepancy for a valuation with center at
$z$.

It is important to remember that the discrepancy $\disc(v) =
\disc(v,X)$ is not birationally invariant. In general, it depends
not only on the valuation $v$, but on the birational model
$X$. An exception to this rule is the following situation.

\begin{defn}\label{crepant}
A dominant, generically one-to-one map $\pi:Y \to X$ is called {\em
crepant} if for every local trivialization $\omega$ of the canonical
bundle $K_X$ the pullback $\pi^*\omega$ is a local trivialization of
the canonical bundle $K_Y$.
\end{defn}

For every proper crepant map $Y \to X$, every valuation $v$ which
centered at a point $\delta(v) \subset X$ also has a center on $Y$,
and we have $\disc(v,X) = \disc(v,Y)$.

\subsection{Galois extensions.}

Assume given a finite Galois extension $K_1/K$ with Galois group
$G$. For every valuation $v$ of the fields $K$, the tensor product
$K_1 \otimes_K K_v$ splits as the sum
\begin{equation}\label{splt}
K_1 \otimes_K K_v = \prod K_{v_i},
\end{equation}
where $K_{v_i}$ are completions of the field $K_1$ with respect to
certain valuations $v_i$. Each of the valuations $v_i$ satisfies
\begin{equation}\label{branch}
v_i(x) = r v(x), \qquad x \in K \subset K_1,
\end{equation}
where $r$ is a positive integer, the same for every $v_i$. This
integer is called the {\em ramification index} of the valuation $v$
in the Galois extension $K_1/K$. If we extend the valuation $v$ to
differential forms, as in \eqref{val.forms}, then the equation
\eqref{branch} holds for every $i$-form $\alpha \in \Omega^i(K/k)
\subset \Omega^i(K_1/k)$.

The valuations $v_i$ are the only valuations of the field $K$ which
satisfy \eqref{branch}. The Galois group $G$ acts transitively on
the set of all the $v_i$. Fix for convenience one of these
valuations, say $v_1$. The stabilizer $G_1 \subset G$ of the
valuation $v_1$ is called the {\em decomposition subgroup} of the
valuation $v$. Each of the completed fields $K_{v_i}$ is a Galois
extension of the field $K_v$ with Galois group conjugate to $G_1
\subset G$. 

One also distinguishes a smaller {\em inertia subgroup} $I \subset
G_0 \subset G$ in the Galois group $G$ which by definition consists
of elements that preserve the valuation $v_1$ {\em and} act
trivially of the residue field $k_{v_1} = O_{v_1}/\m_{v_1}$. Since
all the fields are of characteristic $0$, the inertia subgroup $I
\subset G$ is a cyclic group. Its order is the ramification index
$r$. Moreover, the action on the Zariski cotangent space
$\m_{v_1}/\m_{v_1}^2$ ($1$-dimensional over the residue field
$k_{v_1}$) identifies $I$ with the subgroup $\mu_r \subset
k_{v_1}^*$ of roots of unity in the field $k_{v_1}$. Since the base
field $k$ is algebraically closed, all the roots of unity actually
lie in the subgroup $k^* \subset k_{v_1}$, so that the subgroup $I =
\mu_r$ also lies in the subgroup $k^* \subset k_{v_1}$.

To sum up, every valuation $v$ defines a ramification index $r \in
\Z, r \geq 1$, a decomposition subgroup $G_1 \subset G$ and an
inertia homomorphism $g:\mu_r \to G$. Both $G_1 \subset G$ and
$g:\mu_r \to G$ are defined up to a conjugation within $G$ (this
corresponds to a choice of the valuation $v_1$ in the decomposition
\eqref{splt}).

Since every valuation is multiplicative, we see that the inertia
subgroup $I = g(\mu_r) \subset G$ acts on the $k$-th power
$m_{v_1}^k/\m_{v_1}^{k+1}$ by the $k$-th power of the fundamental
character $\chi:\mu_r \to k^*$. For future use, we will rephrase
this in the following way.

\begin{lemma}\label{chi}
Every eigenvector $x \in K_{v_1}$ of the inertia group $I = g(\mu_r)
\subset G$ satisfies
$$
g(a)(x) = \chi(a)^{v_1(x)}x,
$$
where $a \in \mu_r$ is an arbitrary element in the group $\mu_r$ of
$r$-th roots of unity, and $\chi:\mu_r \to k^*$ is the canonical
embedding.\endproof
\end{lemma}

\section{Monomial valuations}

\subsection{General monomial valuations.}

Let $V$ be a vector space over the base field $k$ of dimension $\dim
V = n$. Assume given an algebraic action of the multiplicative group
$k^*$ on the vector space $V$. Consider the eigenspace decomposition
\begin{equation}\label{splitting}
V = \bigoplus_i V_i
\end{equation}
with respect to the $k^*$-action, and let $a_i \in \Z$ be the
associated weights, so that we have
$$
\lambda \cdot v = \lambda^{-a_i}v, \qquad v \in V_i,\lambda \in k^*.
$$
Every function $f \in O$ in the polynomial ring $O = k[V]$ of
algebraic functions on $V$ has a weight decomposition
$$
f = \sum f_i, \qquad f_i \in K, i \in \Z,
$$
where $f_i$ satisfy $\lambda \cdot f_i = \lambda^i f_i$ for every
$\lambda \in k^*$.

Let $K = K(V)$ be the field of rational functions on $V$. Assume
that the action map $k^* \to GL(V)$ is injective or, equivalently,
that the largest common divisor of the numbers $a_i \in \Z$ equals
$1$. Then for every number $p \in Z$ there exists a rational
function $f \in K$ of weight $p$.

\begin{defn}
The {\em monomial valuation} $v:K^* \to \Z$ of the variety $V$
associated to a $k^*$-action on the vector space $V$ is defined by
\begin{align}\label{monom.eq}
\begin{split}
v(f) &= \min \{i \in \Z | f_i \neq 0\} \qquad f \in O,\\
v(f/g) &= v(f) - v(g) \qquad f,g \in O.
\end{split}
\end{align}
\end{defn}

This is indeed a valuation: the map $v:K^* \to \Z$ is surjective,
and it obviously satisfies the conditions
\thetag{V1}-\thetag{V2}. The residue field $k_v$ for the monomial
valuation $v$ coincides with the field $K(V)^{k^*} \subset K(V)$ of
$k^*$-invariant functions in $K(V)$.

The monomial valuation $v$ has a center $\delta(v) \in V$ if and
only if all the numbers $a_i$ are non-negative, $a_i \geq 0$. Assume
from now on that this is the case.  Say that the monomial valuation
is {\em positive} if the numbers the $a_i$ are strictly positive,
$a_i \geq 1$. Then the center $\delta(v) \in v$ is the closed point
$0 \in V$, and we have $v(f) \geq 1$ if and only if a polynomial $f
\in O$ has no constant term.

Extend the monomial valuation $v$ to differential forms
$\Omega^i(K/k)$ as in \eqref{val.forms}.  The $k^*$-action also
extends to the spaces $\Omega^i(V/k)$ of forms on $V$, and
\eqref{monom.eq} holds for the forms as well as for the
functions. This makes it easy to compute the discrepancy
$\disc(v,V)$. Indeed, the top form $\omega \in \Omega^n(V/k)$ on $V$
given by the determinant is an eigenvector for the $k^*$-action with
eigenvalue
$$
a = \sum_i a_i\dim V_i,
$$
and the discrepancy by definition equals $\disc(v) = a - 1$. The
simplest possible positive monomial valuation is obtained by taking
$V = V_1$ and $a_1 = 1$. This is the adic valuation centered at the
point $0 \in v$. Its discrepancy, just as it should be, is equal to
$n - 1$.

The monomial valuations enjoy the following extremal property.

\begin{lemma}\label{mon.min}
Assume given a valuation $w$ of the variety $V$ with center at $0
\in V$, and assume that for a certain positive monomial valuation
$v$ we have
$$
w(f) \geq v(f)
$$
for every linear function $f \in V^* \subset O$ on $V$. 
\begin{enumerate}
\item The discrepancy satisfies $\disc(w) \geq \disc(v)$.
\item If we have $\disc(w) = \disc(v)$, then the valuations $w$ and
$v$ coincide, $w=v$.
\end{enumerate}
\end{lemma}

\proof{} It is convenient to split each of the subspaces $V_i
\subset V$ in \eqref{splitting} into a sum of one-dimensional
subspaces, so that we have a basis $v_i$ of the vector space $V$, a
collection of integers $a_i$, and for every $i$ we have
$$
\lambda \cdot v_i = \lambda^{-a_i}v_i, \qquad \lambda \in k^*.
$$
Let $x_i$ be the dual basis of the vector space $V^*$ of linear
functions on $V$.

To prove \thetag{i}, note that the differential form $\omega = dx_1
\wedge \cdots \wedge dx_n$ is a trivialization of the canonical
bundle $K_V$ and an eigenvector of the $k^*$-action, with weight
$$
\disc(v) = \sum a_i - 1.
$$
Therefore we indeed have
\begin{align*}
\disc(w) &= w(\omega) - 1 \geq \sum w(dx_i) -1 \geq \sum w(x_i) -1 \\
&\geq \sum v(x_i) - 1 = \sum a_i - 1 = \disc(v).
\end{align*}
To prove \thetag{ii}, it suffices to prove that for every algebraic
function $f \in O$ on $V$ we have $w(f) = v(f)$. First we prove that
$w(f) \geq v(f)$. Indeed, by assumption $w(x_i) \geq v(x_i)$ for all
linear functions $x_i$. Therefore $w(f) \geq v(f)$ for every
monomial $f = x_{i_1} \cdots x_{i_p} \in O$. An arbitrary algebraic
function $f$ is a sum
$$
f = \sum f_p
$$
of monomials $f_p$. Since all monomials are eigenvectors of the
$k^*$-action, we have
$$
v(f) = \min(v(f_p)),
$$
and $w(f) \geq v(f)$. 

To prove the equality, proceed by induction on the degree of the
polynomial $f \in O$. We know that $v(f) = w(f) = 0$ for every
constant $f \in k \subset O$. Assume that $v(f) = w(f)$ for all
polynomials $f \in O$ of degree less than $p$, and let $f$ be a
polynomial of degree $p$. If $f$ has non-zero constant term, then
$v(f) = w(f) = 0$. Therefore we can assume that $f$ has zero
constant term, so that $w(f) \geq v(f) \geq 1$. Note that in this
case $v(f) = v(df)$. We have
$$
df = \sum f_idx_i
$$
for some $f_i \in O$, and
$$
v(df) = \min(v(f_idx_i)) = \min (v(f_i) + v(dx_i)) = \min(v(f_i) +
a_i).
$$
Choose one of the $f_i$ for which the minimum is attained, say,
$f_1$. Since each of $f_i$ is a polynomial of degree less than $p$,
we have $w(f_1) = v(f_1)$ by the inductive assumption. The form
$dx_1 \wedge \cdots \wedge dx_n$ trivializes the canonical bundle
$K_V$, so that we have
\begin{align*}
v(f_1dx_1 \wedge \cdots \wedge dx_n) &= v(f_1) + \disc(v) + 1 =
w(f_1) + \disc(w) + 1 \\
&= w(f_1dx_1 \wedge \cdots \wedge dx_n).
\end{align*}
But since $f_1dx_1 \wedge \cdots \wedge dx_n = df \wedge dx_2 \wedge
\cdots \wedge dx_n$, the right hand side satisfies
\begin{align*}
w(f_1dx_1 \wedge \cdots \wedge dx_n) &\geq w(df) + w(dx_2) + \ldots +
w(dx_n) \\
&\geq w(f) + w(x_2) + \ldots + w(x_n),
\end{align*}
while the left hand side is equal to
\begin{align*}
v(f_1dx_1 \wedge \cdots \wedge dx_n) &= v(f_1) + v(x_1) + \ldots + v(x_n) \\
&= v(f_1dx_1) + v(x_2) + \ldots + v(x_n) \\
&= v(df) + v(x_2) + \ldots + v(x_n) \\
&= v(f) + v(x_2) + \ldots + v(x_n).
\end{align*}
We conclude that
$$
w(f) + w(x_2) + \ldots + w(x_n) \leq v(f) + v(x_2) + \ldots +
v(x_n).
$$
Since we know that $w(x_i) \geq v(x_i)$ for every $i$, this yields
$w(f) \leq v(f)$, which in turn yields $w(f) = v(f)$. This
establishes the induction step and finishes the proof.
\endproof

\begin{corr}\label{adic.minimal}
Assume that a valuation $v$ of an irreducible algebraic variety $X$
has a center $\delta(v) \in X$, and that $X$ is smooth near the
point $\delta(v)$. Let $Z = \overline{\delta(v)} \subset X$ be the
Zariski closure of the point $\delta(v) \in X$. Then the discrepancy
satisfies
$$
\disc(v,X) \geq \codim(Z) - 1,
$$
and the equality holds if and only if $v = v_Z$ coincides with the
$Z$-adic valuation.
\end{corr}

\proof{} Let $V = T_{\delta{v}}X$ be the tangent space. Since the
question is local on $X$, we can replace $X$ with the vector space
$V$, and $\delta(v)$ with $0 \in V$. The $Z$-adic valuation $v_Z$ is
monomial, has discrepancy $\codim(Z) - 1$, and for every linear
function $x$ on $V$ we have $v(x) \geq 1 = v_Z(x)$.
\endproof

If a monomial valuation $v$ is not positive, so that one of the
numbers $a_i$ is equal to $0$, then it is no longer true that $v$ is
centered at $0 \in V$. However, $v$ still has a center $\delta(v)
\in V$. This is the generic point of the subspace $V_0 \subset V$ of
$k^*$-invariant vectors. Replacing the base field $k$ with the field
$K(V_0)$ of rational functions on $V_0$, we immediately generalize
Lemma~\ref{mon.min} to all valuations $w$ of $V$ with center at
$\delta(v) \in V$ and such that $w(f) \geq v(f)$ for every linear
function $f \in V^* \subset O$ on $V$.

\subsection{Actions of a cyclic group.}\label{cyclic}

Let now $V$ be an $n$-dimensional vector space over the field $k$,
let $r \in \Z, r \geq 2$ be an integer, and let $\mu_r \subset k^*$
be the cyclic group of roots of unity of order $r$. Assume given an
injective homomorphism $g:\mu_r \hookrightarrow GL(V)$. The
homomorphism $g$ canonically extends to a $k^*$-action on $V$ in the
following way. Let
$$
V = \bigoplus V_i
$$
be the eigenspace decomposition of the vector space $V$, so that
$\mu_r$ act on $V_i \subset V$ by the $(-a_i)$-th power $\chi^{-a_i}$
of the fundamental character $\chi:\mu_r \to k^*$. Choose the
integers $a_i$ so that $0 \leq a_i < r$, and let
$$
\lambda \cdot v = \lambda^{-a_i}v, \qquad \lambda \in k^*, v \in
V_i.
$$
Since $g:\mu_r \hookrightarrow GL(V)$ is injective, the largest
common denominator of the numbers $a_i$ is equal to $1$. Therefore
we obtain a monomial valuation $v$ of the field $K = K(V)$, which we
denote by $v = v_g$. For every polynomial $f \in O = k[V^*]$ which
is an eigenvector of the group $\mu_r$, the valuation $v_g$
satisfies 
$$
g(x) \cdot f = \chi(x)^{v_g(f)}f, \qquad x \in \mu_r.
$$
Let $X = V/\mu_r$ be the quotient of the vector space $V$ by the
$\mu_r$ action, and let $K_0 = K^{\mu_r}$ be the field of rational
functions on $X$. By restriction, we obtain a map $v_g:K_0^* \to
\Z$. This map is a surjection onto $r\Z \subset \Z$. Therefore the
quotient $v_g/r$ is a well-defined valuation of the field $K_0$.

If the cyclic group $g(\mu_r) \subset SL(V) \subset GL(V)$ preserves
the determinant, then the quotient $X = V/\mu_r$ has Gorenstein
singularities. Indeed, the canonical bundle $K_X$ is trivial. The
determinant form $\omega \in \Omega^n(V/k)$ is $\mu_r$-invariant and
descends to a trivialization of the bundle $K_X$. In this case it
makes sense to speak of the discrepancy of the valuation $v_g/r$ of
the quotient variety $X$. By definition, it equals
\begin{equation}\label{disc.quot}
\disc\left(\frac{v_g}{r},X\right) = \frac{v_g}{r}(\omega) - 1 =
\frac{1}{r}(\disc(v_g,V)+1) -1 = \sum \frac{a_i}{r}\dim V_i -1.
\end{equation}
The number
$$
\age(g) = \sum \frac{a_i}{r}\dim V_i
$$
is called the {\em age} of the homomorphism $g:\mu_r \to SL(V)$. The
age is a positive integer, $\age{g} \in \Z$, $\age(g) \geq 1$.

\subsection{Symplectic actions of an arbitrary finite group.}

More generally, assume that a finite group $G$ acts on the vector
space $V$, the map $G \to GL(V)$ is injective, and we are given an
injective homomorphism $g:\mu_r \hookrightarrow G$.  By restriction,
we obtain a $\mu_r$-action on the vector space $V$. As in the last
subsection, denote the associated monomial valuation of $V$ by
$v_g$.

Consider the quotient $X = V/G$ and the field of rational functions
$K_0 = K(X) = K^G$. As in the case $\mu_r = G$, we would like to
define a valuation of the variety $X$ associated to the monomial
valuation $v_g$. By restriction, we obtain a map $K_0^* \to \Z$,
which is a surjection on a subgroup $r_g\Z \subset \Z$, where $r_g$
is a certain integer. The quotient $v_0 = v_g/r_g$ is a well-defined
valuation of the field $K_0 \subset K$. 

The extension $K/K_0$ is Galois with the Galois group $G$. Since
$v_g$ extends the valuation $v_0$ to the field $K$, the ramification
index of the extension $K/K_0$ is equal to $r_g$, and we have the
inertia homomorphism $I:\mu_{r_g} \to G$. 

\begin{lemma}
The subgroup $g(\mu_r) \subset G$ lies in the inertia subgroup
$I(\mu_{r_g}) \subset G$. Thus the ramification index $r_g$ must be
a multiple of the integer $r$.
\end{lemma}

\proof{} By definition of the monomial valuation $v = v_g$, the
subgroup $g(\mu_r) \subset G$ preserves the valuation $v:K(V)^* \to
\Z$. Therefore $g(\mu_r) \subset G_1$ lies in the decomposition
subgroup $G_1 \subset G$. To prove that $g(\mu_r) \subset
I(\mu_{r_g})$ lies in the inertia subgroup $I(\mu_{r_g}) \subset
G_1$, we have to prove that $g(\mu_r)$ acts trivially on the residue
field $k_v$. 

But the residue field coincides with the subfield $K(V)^{k^*}
\subset K(V)$ of $k^*$-invariant rational function. By definition of
the monomial valuation, this $k^*$-action on $V$ extends the action
of the subgroup $\mu_r \subset k^*$. Therefore every $k^*$-invariant
function is also $g(\mu_r)$-invariant.
\endproof

In general, it might happen that $r_g \neq r$.  We show that if $G
\in Sp(V)$ preserves a symplectic form on $V$ (the case that we
study in this paper), this never happens.

\begin{lemma}\label{r=r.g}
If $V \cong V^*$ as representations of the group $G$, then the
ramification index $r_g$ coincides with the integer $r$.
\end{lemma}

\proof{} Since the inertia group $I(\mu_{r_g}) \subset G$ is
commutative, and $g(\mu_r) \subset I(\mu_{r_g})$, the eigenspace
decomposition
$$
V = \bigoplus V_i
$$
with respect to $\mu_r$ is invariant under $I(\mu_{r_g})$. Therefore each
of the spaces $V_i$ is generated by eigenvectors of the
$I(\mu_{r_g})$-action. 

Let $v \in V_i$ be such an eigenvector. Applying Lemma~\ref{chi} to
linear functions on $V$, we see that for any element $a \in
\mu_{r_g}$ we have
$$
I(a) \cdot v = \chi(a)^{-a_i}v.
$$
This number $a_i$ does not depend on the choice of the eigenvector
$v \in V_i$.  Since $V_i$ is generated by these eigenvectors, the
group $I(\mu_{r_g})$ acts on every vector $v \in V_i$ by the same
character $\chi^{-a_i}$.

We see that all the non-zero weights $a_i$ of the
$I(\mu_{r_g})$-action on the vector space $V$ lie in the interval
$(0,r)$. Therefore the non-zero weights of the dual action on $V^*$
lie in the interval $(r_g-r,r)$. But by assumption $V \cong
V^*$. This is possible only if $r_g = r$.
\endproof

Since the symplectic group $Sp(V) \subset SL(V)$ preserves the
determinant, for every finite subgroup $G \subset Sp(V)$ the
quotient $X = V/G$ has Gorenstein singularities. As in
\eqref{disc.quot}, the discrepancy of the monomial valuation $v_g/r$
is equal to
$$
\disc\left(\frac{v_g}{r},X\right) = \age(g) - 1.
$$
In our situation $G \subset Sp(V)$ this number is very easy to
compute. 

\begin{lemma}\label{age.onehalf}
Assume that $V \cong V^*$ as representations of the group $G$.  The
age $\age{g}$ of an arbitrary homomorphism $g:\mu_r \to G$ equals
$$
\age{g} = \frac{1}{2}\codim V_0,
$$
where $V_0 = V^g \subset V$ is the subspaces of vectors invariant
under $g(\mu_r) \subset G$.
\end{lemma}

\proof{} Let $a_i$ be the full set of weights of the $\mu_r$-action
on $V$ induced by $g:\mu_r \to G$, so that $a_0 = 0$, $0 < a_i < r$
for $i > 0$ and
$$
\age(g) = \frac{1}{r}\sum a_i\dim V_i.
$$
Since $V \cong V^*$, the multiplicity $\dim V_i$ of every non-zero
weight $a_i$ coincides with the multiplicity of the weight $r -
a_i$. Therefore
\begin{align*}
2\age{g} &= \frac{1}{r} \sum_{i > 0}a_i\dim V_i + \frac{1}{r} \sum_{i
> 0} (r-a_i)\dim V_i = \frac{1}{r} \sum_{i > 0}(a_i + (r-a_i))\dim V_i\\ 
&= \sum_{i > 0}\dim V_i = \codim V_0.\hspace{0.55\linewidth}
\square
\end{align*}

\noindent
{\bf Remark.} We note that in the general situation $G \subset SL(V)$ both
Lemma~\ref{r=r.g} and Lemma~\ref{age.onehalf} are false. The latter
is not surprising. Indeed, in the language of McKay correspondence
developed in \cite{R}, Lemma~\ref{age.onehalf} corresponds to the
fact that any smooth crepant resolution $Y \to X = V/G$ is
semismall. This is usually false unless $G \subset Sp(V)$. 

However, the failure of Lemma~\ref{r=r.g} is rather alarming. The
simplest example of this failure is the case when $V = \C^4$, and $G
= \mu_4 \cong \Z/4\Z$ acts on each $\C \subset V$ by the fundamental
character $\chi$. If one takes the natural embedding $g:\mu_2 \to
\mu_4 = G$, then $r_g = 4$ and $r = 2$.

Fortunately, this is also the simplest case when the quotient $V/G$
admits no smooth crepant resolution (see \cite{IN}). One can
conjecture that Lemma~\ref{r=r.g} holds whenever a smooth crepant
resolution $Y \to X = V/G$ does exist. This can be formulated in the
following purely algebraic way.

\begin{conjecture}
Assume that the quotient $X = V/G$ of a vector space $V$ by a finite
subgroup $G \subset SL(V)$ admits a smooth crepant resolution $Y \to
X$.

Then for every cyclic subgroup $g(\mu_r) \subset G$, there exists a
$G$-invariant rational function $f \in K(V)$ on $V$ which has weight
$1$ with respect to the induced $g(\mu_r)$-action. In other words,
we have
$$
g(a) \cdot f = \chi(a)f
$$
for every element $a \in \mu_r$.
\end{conjecture}

\subsection{Reid's conjecture.}

We are now in a position to give the precise formulation of the
conjecture from the paper \cite{R} which we will prove.

Let $V$ be a vector space, and let $G \subset SL(V)$ be a finite
subgroup. Consider the quotient $X = V/G$. Assume given a proper
smooth crepant resolution $\pi:Y \to X$.

Since the quotient $X = V/G$ and the resolution $Y$ are not compact,
it is convenient to work with the Borel-Moore homology groups
$H^c_\idot(Y,\Z)$ instead of the ordinary holomogy groups
$H_\idot(Y,\Q)$. These are by definition the relative homology
groups
$$
H^c_\idot(Y,\Q) = H_\idot\left(\wt{Y},\pt;\Q\right)
$$
of the one-point compactification $\wt{Y}$ modulo the infinite point
$\pt \in \wt{Y}$. The Borel-Moore homology is dual to the cohomology
with compact support. The advantage over the ordinary homology is
the following: every algebraic cycle $Z \subset Y$ has a
well-defined fundamental class $\cl(Z) \in H^c_\idot(Y,\Q)$.

For every homomorphism $g:\mu_r \to G$, consider the monomial
valuation $v_g/r_g$ of the quotient variety $X$. Since $Y \to x$ is
proper, the valuation $v_g/r_g$ has a center $\delta(v_g/r_g) \in
Y$. Denote by $Z_g = \overline{\delta(v_r/r_g)} \subset Y$ the
Zariski closure of this center. Note that for two homomorphisms
$g_1:\mu_{r_1} \to G$, $g_2:\mu_{r_2} \to G$ we have $Z_{g_1} =
Z_{g_2}$ if and only if $r_1 = r_2$ and the homomorphisms
$g_1,g_2:\mu_{r_1} \to G$ are conjugate by an element of $G$.

\begin{conjecture}\label{main}
The fundamental classes $\cl(Z_g)$ of the algebraic cycles $Z_g
\subset Y$ form a basis of the homology $\Q$-vector space
$H^c_\idot(Y,\Q)$.  Moreover, for every $g:\mu_r \to G$, we have
$\age(g) = \codim(Z_g)$.
\end{conjecture}

This formulation makes sense for an arbitrary subgroup $G \subset
SL(V)$, but we shall prove it only for $G \subset Sp(V)$. Whether
this is the correct statement for the general situation $G \subset
SL(V)$, we do not know.

\section{Crepant resolutions}

\subsection{Topology of crepant resolutions.}

Let $V$ be a finite-dimensional vector space equipped with an action
of a finite group $G \subset V$. Let $X = V/G$ be the quotient
variety. Assume given a smooth proper crepant resolution $\pi:Y \to
X$.

\label{strt}
The vector space $V$ is naturally stratified by subspaces $V^H
\subset V$ of $H$-invariant vectors for various subgroups $H \subset
G$ (some strata may coincide, since we might have $V^{H_1} = V^{H_2}
\subset V$ for different subgroups $H_1,H_2 \subset G$).  This
stratification induces a stratification $X_H$ of the quotient
variety $X = V/G$. For each subgroup $H \subset G$, the complement
$$
X_H^o = X_H \setminus \bigcup_{X_{H_1} \varsubsetneq X_H} X_{H_1}
$$
is a smooth locally closed subvariety $X_H^o \subset X$. For an
arbitrary cyclic group homomorphism $g:\mu_r \to G$, we will denote
by $X_g \subset X$ the stratum corresponding to the subgroup
$g(\mu_r) \subset G$.

Let $Y_H = \pi^{-1}(X_H) \subset Y$ be the associated stratification
of the resolution $Y$. We recall the following general property of
smooth crepant resolutions of symplectic quotient singularities.

\begin{prop}[\refsemismall]\label{semismall}
For every stratum $X_H \subset X$, the map $\pi:Y_H \to X_H$ is a
locally trivial fibration in \'etale topology over the open dense
subset $X_H^o \subset X_H$, and we have
$$
\codim(Y_H) \geq \frac{1}{2}\codim(X_H).
$$
In particular, the resolution $\pi:Y \to X$ is semismall for the
stratification $X_H$. \endproof
\end{prop}

\begin{remark}
The statement of \refsemismall{} uses a slightly different
description of the stratification on $X$ which does not mention
explicitly the subgroups $H \subset G$. Nevertheless, it is very
easy to check that this is the same stratification.
\end{remark}

We will need the following corollary of Proposition~\ref{semismall}.

\begin{defn}\label{max.cycle}
An irreducible closed subvariety $Z \subset Y$ is called a {\em
maximal cycle} if 
$$
\codim(Z) = \frac{1}{2}\codim(\pi(Z)).
$$
\end{defn}

\begin{corr}\label{ss}
Every maximal cycle $Z \subset Y$ is an irreducible component of a
certain stratum $Y_H \subset Y$.
\end{corr}

\proof{} Let $X_H^o \subset X$ be the stratum which contains the
generic point of the closed subvariety $\pi(Z) \subset X$. Since $Z$
is irreducible, we have $Z \subset Y_H \subset Y$. Therefore the
dimension $\dim(Z/\pi(Z))$ of the generic fiber of the projection
$\pi:Z \to \pi(Z)$ satisfies
$$
\dim (Z/\pi(Z)) \leq \dim (Y_H /X_H ) \leq \frac{1}{2}\codim(X_H).
$$
Since $\pi(Z) \subset X_H$, we have $\codim(\pi(Z)) \geq
\codim(X_H)$ and
\begin{align*}
\codim(Z) &= \codim(\pi(Z)) - \dim(Z/\pi(Z)) \geq \codim(\pi(Z)) -
\frac{1}{2}\codim(X_H) \\
&\geq \codim(\pi(Z)) - \frac{1}{2}\codim(\pi(Z))
= \frac{1}{2}\codim(\pi(Z)).
\end{align*}
Sine the cycle $Z \subset Y$ is maximal, all the inequalities are in
fact equalities. Therefore $\pi(Z) = X_H \subset X$, and $Z \subset
Y_H$ is an irreducible component of maximal possible dimension.
\endproof

Our proof of Conjecture~\ref{main} in the symplectic case splits
into three steps.
\begin{enumerate}
\renewcommand{\labelenumi}{\arabic{enumi}.}
\item We prove that for every $g:\mu_r \to G$, the subvariety $Z_g
\subset Y$ is a maximal cycle.
\item We prove that every maximal cycle $Z \subset Y$ coincides with
$Z_g$ for some homomorphism $g:\mu_r \to G$.
\item We prove that the homology group $H^c_\idot(Y,\Q)$ is a
$\Q$-vector space generated by the classes of the maximal cycles
$Z_g$.
\end{enumerate}

Steps 1 and 2 are the subject of next two subsections. Step 3 is
purely topological, and we deal with it in the next section.

\subsection{From monomial valuations to maximal cycles.}

Let $g:\mu_r \to G$ be the embedding of a cyclic subgroup $\mu_r
\subset g$. Consider the associated monomial valuation $v_g/r$ of
the quotient variety $X$. The valuation $v_g/r$ has a center
$\delta(v_g/r,X) \in X$ in $X$, which coincides with the generic
point of the stratum $X_g \subset X$.

Let $Z_g \subset Y$ be the Zariski closure of the center
$\delta(v_g/r,Y) \subset Y$ of the valuation $v_G/r$ in $Y$. Since
$\pi(\delta(v_g/r,Y)) = \delta(v_g/r,X)$, we see that $\pi(Z_g) = X_g
\subset X$.

\begin{prop}\label{mon=>max}
The subvariety $Z_g \subset Y$ is a maximal cycle of codimension
$\codim Z_g = \age(g)$. Moreover, the valuation $v_g/r$ of the
variety $Y$ coincides with the $Z_g$-adic valuation.
\end{prop}

\proof{} Let $v = v_g/r$. Since $\pi:Y \to X$ is crepant, the
discrepancy $\disc(v) = \disc(v,X) = \disc(v,Y)$ is the same for $X$
and for $Y$. By Lemma~\ref{age.onehalf} it equals
$$
\disc(v) = \age(g) - 1 = \frac{1}{2}\codim X_g - 1.
$$
On the other hand, by Corollary~\ref{adic.minimal} we have
$$
\disc(v) \geq \codim(Z_g) - 1.
$$
Therefore $2\codim(Z_g) \leq \codim X_g$. Applying
Proposition~\ref{semismall}, we conclude that
$$
\codim(Z_g) = \frac{1}{2}\codim(X_g).
$$
Thus $Z_g \subset Y$ is a maximal cycle. Moreover, we have $\disc(v)
= \codim(Z_g) - 1$, and $v$ coincides with the $Z_g$-adic valuation
by the second claim of Corollary~\ref{adic.minimal}.
\endproof

\subsection{From maximal cycles to monomial valuations.}

Assume given a maximal cycle $Z \subset Y$ with generic point $z
\subset Y$. The closure $X_Z \subset \overline{\pi(z)} \subset X$ is
a certain stratum in the variety $X = V/G$. Denote by $v = v_Z$ the
$Z$-adic valuation of the field $K(X) = K(Y)$. Since the resolution
$Y$ is smooth and crepant, we have
$$
\disc(v,X) = \disc(v,Y) = \codim(Z) - 1 = \frac{1}{2}\codim X_z - 1.
$$
Let $r \in \Z$ be the ramification index of the Galois extension
$K(V)/K(X)$, and let $g:\mu_r \to G$ be the inertia
homomorphism. Consider the decomposition
$$
V = \bigoplus_{i \geq 0}V_i
$$
of the vector space $V$ with respect to the induced $\mu_r$-action,
and let
$$
V^* = \bigoplus_{i \geq 0}V_i^*
$$ 
be the dual decomposition of the space $V^*$ of linear functions on
$V$. Let $a_i$ be the weight of the subspace $V_i \subset V$, $a_0 =
0$, $0 < a_i < r$ for $i \geq 1$. Consider the monomial valuation
$v_g$ of $V$ associated to the homomorphism $g:\mu_r \to G$.

Extends the valuation $v$ to a valuation $v_0$ of the variety $V$ so
that $v_0 = rv$ on $K(X) \subset K(V)$. The discrepancy $\disc(v_0)$
equals
$$
\disc(v_0) = r(\disc(v)+1)-1.
$$
Moreover, by Lemma~\ref{chi} we have
$$
v_0(x) = v_g(x) = a_i \mod(r)
$$
for every linear function $x \in V_i^*$. Since $v_0(x) \geq 0$, this
implies that $v_0(x) \geq v_g(x)$. By the definition of the monomial
valuation, we must have
\begin{equation}\label{geq}
v_0(x) \geq v_g(x)
\end{equation}
for an arbitrary linear function $x \in V^*$.

\begin{prop}\label{max=>mon}
The $Z$-adic valuation $v = v_Z$ of the variety $Y$ coincides with
the monomial valuation $v_g/r$.
\end{prop}

\proof{} Let $V_Z \subset V$ be the subspace associated to the
stratum $X_Z = \pi(Z) \subset X$. We begin by proving that $V_Z$
coincides with the subspace $V_0 = V^g \subset V$ of
$g(\mu_r)$-invariant vectors in $V$.

The generic point of the subspace $V_Z \subset V$ is the center
$\delta(v_0)$ of the valuation $v_0$ of the variety $V$. By
definition the inertia subgroup acts trivially on the residue field
of the valuation $v_0$. Therefore $V_Z \subset V_0$.

Choose a complement $V' \subset V_0$ to $V_Z \subset V_0$. Consider
the decomposition
$$
V = V_Z \oplus V' \oplus \bigoplus_{i \geq 1}V_i,
$$
and define a monomial valuation $v'$ by assigning weight $0$ to
$V_Z \subset V$, weight $r$ to $V' \subset V$ and weight $a_i$ to
$V_i$, $i \geq 1$. 

By Lemma~\ref{chi}, for every linear function $x \in V^*$ which is
an eigenvector of the $\mu_r$-action we have
$$
v_0(x) = v'(x) \mod(r)
$$
Moreover, if a linear function $x \in V^*$ vanishes on $V_Z \subset
V$, we have $v_0(x) \geq 1$. Since $v_0(x) = 0 \mod(r)$, this
implies $v_0(x) \geq r$. Therefore for every linear function $x \in
V^*$ we have
$$
v_0(x) \geq v'(x).
$$
By Lemma~\ref{mon.min} this implies that
$$
\disc(v_0) \geq \disc(v').
$$
But
$$
\disc(v_0) = r(\disc(v)+1) - 1 = \frac{r}{2}\codim V_Z - 1,
$$
while by Lemma~\ref{age.onehalf}
$$
\disc(v') = \frac{r}{2}\codim V_0 + r\dim V' - 1 = \frac{r}{2}\codim
V_Z + \frac{r}{2}\codim V' - 1.
$$
Therefore $\dim V' = 0$, and we indeed have $V_Z = V_0$. 

Moreover, $\disc(v_0) = \disc(v)$. Since $V'=0$, the valuation $v' =
v_g$ coincides with the monomial valuation associated to the inertia
homomorphism $g:\mu_r \to G$. Applying
Lemma~\ref{mon.min}~\thetag{ii}, we conclude that $v_0 = v' = v_g$
and $v = v_0/r = v_g/r$.
\endproof

As a corollary, we see that the stratum $X_Z = \pi(Z_g) \subset X$
coincides with the stratum $X_g \subset X$ associated to the
homomorphism $g:\mu_r \to G$.

\section{Homology of a crepant resolution}

Let $V$, $G$, $X$ and $Y$ be as in the last section. To finish the
proof of Conjecture~\ref{main}, it remains to do the third and last
step, namely, to prove the following statement.

\begin{prop}\label{third}
The homology group $H^c_\idot(Y,\Q)$ is a $\Q$-vector space
generated by classes of maximal cycles $Z \subset Y$.
\end{prop}

First we list the necessary topological properties of the quotient
variety $X = V/G$. Recall (see page \pageref{strt}) that the variety
$X$ is naturally stratified by closed subvarieties $X_H \subset X$,
numbered by the subgroups $H \subset G$. The stratum $X_H =
\eta(V^H) \subset X$ is the image of the subspace $V^H \subset V$ of
$H$-in\-va\-ri\-ant vectors under the quotient map $\eta:V \to X$.

\begin{lemma}\label{quot}\mbox{}
\begin{enumerate}
\item Every perverse sheaf of $\Q$-vector spaces on $X$ which is
smooth along all the open strata $X_H^o \subset X$ is a direct sum
of Goersky-MacPherson sheaves supported on strata $X_H \subset X$.
\item For every Goresky-MacPherson perverse sheaf $\K$ supported on a
stratum $X_H \subset X$, the shifted complex $\K[-\dim X_H]$ is an
ordinary sheaf on $X$.
\item We have $H_c^i(X_H,\K) = 0$ unless $i = \dim X_H$.
\end{enumerate}
\end{lemma}

\proof{} Since the vector space $V$ is symplectic, every stratum
$X_H = \eta(V^H) \subset X$ is even-dimensional. This implies
\thetag{i}.

Let $H \subset G$ be an arbitrary subgroup, and let $V^H \subset V$
be the subspace of $H$-invariant vectors in $V$. Denote by $G_0 =
\Stab(V^H)/\Cent(V^H)$ the quotient of the subgroup $\Stab(V^H)
\subset G$ of element $g \in G$ which preverve the subspace $V^H
\subset V$ by the subgroup $\Cent(V^H) \subset \Stab(V^H)$ of
elements $g \in \Stab(V^H)$ which act as identity on $V^H \subset
V$. Then the stratum $X_H \subset X$ is in fact the quotient of the
subspace $V^H \subset V$ by the group $G_0$. Since all the smaller
strata $X_{H'} \subset X_H$ are of codimension $\geq 2$, this
implies that the fundamental group of the open smooth part $X_H^o
\subset X_H$ coincides with the group $G_0$. Therefore every
Goresky-MacPherson sheaf $\K$ supported on $X_H$ is a direct summand
of the sheaf $\eta_{H*}\Q_H$, where $\Q_H$ is the constant perverse
sheaf on $V^H$, and $\eta_H:V^H \to X_H$ is the quotient map.

Thus it suffices to prove \thetag{ii} and \thetag{iii} for the
perverse sheaf $\K = \eta_*\Q_H$. Equivalently, it suffices to prove
both statements for the sheaf $\Q_H$ itself, with $X$ replaced by
$V$ and $X_H$ replaced by $V^H$. In this setting \thetag{ii} and
\thetag{iii} are obvious.  \endproof

We can now begin the proof of Proposition~\ref{third}.

\begin{lemma}
Let $n = \dim X$. For every $k \geq 0$ we have
$$
H_c^{n+k}(Y,\Q) = \bigoplus_{X_H} H_c^{2k}(X_H,R^{n-k}\pi_*\Q),
$$
where the sum is taken over all the strata $X_H \subset X$ with
$\dim X_H = k$. (In particular, $H_c^p(X,\Q) = 0$ for $p < n$.)
\end{lemma}

\proof{} Since the resolution $\pi:Y \to X$ is semismall and locally
trivial over the open strata $X_H^o \subset X$, the direct image
$$
R^\highdot\pi_*\Q[n]
$$
of the constant perverse cheaf $\Q[n]$ on $Y$ is a perverse sheaf on
$X$ smooth along the open strata $X_H^o \subset X$. By
Lemma~\ref{quot}~\thetag{i} (alternatively, by the Decomposition
Theorem), we have a direct sum decomposition
$$
R^\highdot\pi_*\Q = \bigoplus_H \K_H[\dim X_H - \dim X],
$$
where $\K_H$ is a certain Goresky-MacPherson perverse sheaf
supported on the stratum $X_H \subset X$. By
Lemma~\ref{quot}~\thetag{ii}, this implies that
$$
R^\highdot\pi_*\Q = \bigoplus_k R^k\pi_*\Q
$$
in the derived category of complexes with constructible cohomology
on $X$. Therefore the Leray spectral sequence for the map $\pi:Y \to
X$ degenerates, and we have
$$
H_c^p(Y,\Q) = \bigoplus_k H_c^{p-k}(X,R^k\pi_*\Q).
$$
Moreover, we have
$$
R^k\pi_*\Q = \bigoplus_{X_H} \K_H[k-n],
$$
where the sum is taken over all the strata $X_H \in X$ with $\dim
X_H = n-k$. By Lemma~\ref{quot}~\thetag{iii}, this implies that
$$
H_c^{p-k}(X,R^k\pi_*\Q) = \bigoplus_{\dim X_H = n-k}
H_c^{p-k}(X_H,\K_H[k-n])
$$
vanishes unless $p-k = 2(n-k)$, which yields the claim of the lemma.
\endproof

To finish the proof of Proposition~\ref{third}, it suffices to prove
the following.

\begin{lemma}
For every stratum $X_H \subset X$ of dimension $\dim X_H = k$, the
subspace
$$
H_c^{2k}(X_H, R^{n-k}\pi_*\Q) \subset H_c^{n+k}(Y,\Q)
$$
is freely generated by classes of maximal cycles $Z \subset Y$ which
dominate $X_H \subset X$.
\end{lemma}

\proof{} Consider the preimage $Y_H = \pi^{-1}(X_H) \subset Y$. By
proper base change, we can replace the map $\pi:Y \to X$ with the
restriction $\pi:Y_H \to X$. We begin by proving that
\begin{equation}\label{tp}
H_c^i(X_H, R^j\pi_*\Q) = 0
\end{equation}
unless either $j = n-k$, or $i + j < n+k-1$. 

Indeed, for every $j > n-k$ the sheaf $R^j\pi_*\Q$ is supported on
the strata $X_{H'} \subset X_H$ of dimension $\dim X_{H'} = n-j <
k$.  Lemma~\ref{quot}~\thetag{iii} yields \eqref{tp} unless $i =
2(n-j)$ and $i+j = 2n-j = n+\dim X_{H'}$. Since all the strata are
even-dimensional, this implies $i+j < n+k - 1$.

On the other hand, if $j < n-k$, then the sheaf $R^j\pi_*\Q$
vanishes unless $j$ is even, which forces $j < n-k-1$. But $\dim X_H
= k$, so that we have \eqref{tp} unless $i \leq 2k$. Together this
again yields $i+j < n+k - 1$.

We conclude that the Leray spectral sequence for the map $\pi:Y_H
\to X_H$ degenerates for $H_c^{n+k}$, and we have
$$
H_c^{2k}(X_H, R^{n-k}\pi_*\Q) = H_c^{n+k}(Y_H,\Q).
$$
Since $2\dim Y_H \leq n+k$, this group is generated by cohomology
classes of those irreducible components $Z \subset Y_H$ of the
variety $Y_H$ for which we have the equality $2\dim Z = n+k$. These
are precisely the maximal cycles in $Y$ dominating $X_H \subset X$.
\endproof

This proves Proposition~\ref{third} and finishes the proof of
Conjecture~\ref{main} in the symplectic case.

\subsection*{Acknowledgements.} I am grateful to T. Bridgeland, A. Kuznetsov,
M. Leenson and M. Reid for interesting and valuable discussions. I
would like to specifically thank M. Reid for attracting my attention
to the paper \cite{IR}. I would like to thank G. Zigler for
attacting my attention to the paper \cite{DHZ}. I am especially
grateful to M. Verbitsky for many discussions and for many important
suggestions. In particular, it was he who explained to me that,
because of the semi-smallness of smooth crepant resolutions, the
general problem of the McKay correspondence might be solvable in the
case of symplectic quotient singularities.

I am also grateful to V. Ginzburg, who found a serious mistake in an
earlier version of the proof of Proposition~\ref{third}, and to
M. Kontsevich, who brought to my attention important papers
\cite{Ba}, \cite{DL} and corrected my misapprehensions on the
current state-of-the-art in McKay correspondence.

\end{document}